\documentclass{dcds-B}

\def\currentvolume{3}
\def\currentissue{2}
\def\currentyear{2003}
\def\currentmonth{May}
\def\ppages{255--262}

\newtheorem{lemma}{Lemma}[section]
\newtheorem{theorem}[lemma]{Theorem}
\newtheorem{remark}[lemma]{Remark}

\newtheorem{example}[lemma]{Example}

\renewcommand{\subjclassname}{\textup{2000} Mathematics Subject
 Classification}

\renewcommand{\labelenumi}{(\roman{enumi})}
\title[Recurrent Motions in the Navier-Stokes System]
{Recurrent Motions  in the   Nonautonomous Navier-Stokes System}

\author[V. P. Bongolan-Walsh]{Vena Pearl Bongolan-Walsh}
\address[V. P. Bongolan-Walsh]{%
Department of Applied Mathematics \\ Illinois Institute of
Technology\\ Chicago, IL 60616, USA}
 \email[V. P. Bongolan-Walsh]{bongven@iit.edu}

\author[D. Cheban]{David ~Cheban}
\address[D. Cheban]{%
State University of Moldova\\ Department of Mathematics and
Informatics\\ A. Mateevich Street 60\\ MD--2009 Chi\c{s}in\u{a}u,
Moldova} \email[D. Cheban]{cheban@usm.md}

\author[J. Duan]{Jinqiao ~Duan}
\address[J. Duan]{%
Department of Applied Mathematics \\ Illinois Institute of
Technology\\ Chicago, IL 60616, USA}
\address[J. Duan]{%
Department of Mathematics \\ University of Science and Technology of China\\
 Hefei, Anhui 230026 CHINA} \email[J.~Duan]{duan@iit.edu}

\subjclass{Primary:34C35, 34D20, 58F10, 58F12;
Secondary: 35B35, 35B40.}
\keywords{Nonautonomous dynamical
system, skew--product flow,     Navier-Stokes equations,
recurrent motion, Poisson stable motion.}

\begin{document}

 \maketitle

\medskip

\centerline{(Communicated by P. E. Kloeden)}

 \medskip

\begin{abstract}
We prove the existence of recurrent or  Poisson stable
motions in the   Navier-Stokes
 fluid system under  recurrent  or  Poisson stable forcing, respectively.
We use an approach
based on nonautonomous dynamical systems ideas.

\end{abstract}

\section{Introduction}

The mankind has  been  fascinated by  time-periodic,
quasi-periodic, almost periodic and recurrent motions  for
centuries. These motions have been observed in the solar system
(e.g., the Earth rotates around the Sun), and in   fluid systems,
 (e.g., vortices), and in other natural systems. In fact, in   the atmospheric and
oceanic flows, these motions (such as hurricanes) dominate some
weather or climate evolution. From a dynamical system point of
view, these motions are special invariant sets which provide a
part of the skeleton for understanding the global dynamics.
Therefore, the study of time-periodic, quasi-periodic, almost
periodic and recurrent motions is not only   physically useful but
also   mathematically interesting.

However,   rigorous  justification of the existence of such
motions in fluid systems is a much recent issue.   For example,
time-periodic solutions for the two-dimensional (2D) Navier-Stokes
fluid system under external time-periodic forcing  are shown to
exist by  Foias and Prodi \cite{Prodi, Foias}, and Yudovic
\cite{Yu} in 1960s. Note that a time-periodic  external forcing
does not  necessarily induce a time-periodic motion in a dynamical
system, as in the resonance case of  a periodically forced
harmonic oscillator which only has unbounded, non-periodic
motions.

For sufficiently small external nonautonomous forcing  and/or for
sufficiently large  viscosity,   Chepyzhov and Vishik \cite{Vishik} have
investigated the existence of  quasi-periodic and almost periodic
solutions of the 2D Navier-Stokes fluid system.
Duan and Kloeden \cite{Duan, DuanKloeden} have studied periodic,
 quasi-periodic, and almost periodic motions for
large scale quasigeostrophic fluid flows under nonautonomous forcing.

In this paper, we investigate the existence of  recurrent  and
Poisson stable motions (see definitions in the next section)  of
2D Navier-Stokes fluid system under external     forcing. We use a
nonautonomous dynamical systems approach.

The problem of existence  of recurrent and Poisson stable motions for ordinary
differential equations was studied by Shcherbakov \cite{Scer1, Scer2}.

\section{Nonautonomous dynamical systems }

We first introduce some basic concepts in autonomous dynamical
systems.

Let $\mathbb T = \mathbb R $ or $\mathbb R_{+}$,  $X$ be a metric space with metric $\rho$,
and $(X,\mathbb
T,\pi)$ be an autonomous dynamical system. That is, the continuous
mapping $\pi :\mathbb T \times X \to  X$ satisfies the identity
property $\pi (0, x)=x$ and the flow or semi-flow property
 $\pi (t+\tau, x)= \pi(t, \pi(\tau, x))$. The solution mapping for
 the initial-value problem of an autonomous ordinary differential equation
 usually satisfies these properties.

The point $u\in X$ is called a stationary ($\tau$-periodic, $\tau
>0, \tau \in \mathbb T$) point, if $ut=u$ ($u\tau = u$
respectively) for all $t\in \mathbb T$, where $ut:=\pi (t,u)$.

The number $\tau \in \mathbb T$ is called an $\varepsilon >0$ shift
(almost period) of point $u \in X$ if $\rho (u\tau,u)<\varepsilon
$ (respectively $\rho (u(\tau +t),ut)<\varepsilon$ for all $t\in \mathbb
T$).

The point $u \in X $ is called almost recurrent (almost periodic,
in the sense of Bohr) if for any $\varepsilon$ there exists a
positive number $l$ such that on any segment of length $l$, there
is an $\varepsilon$ shift (almost period) of point $u\in X$.

If a point $u\in X$ is almost recurrent and the hull
$H(u)=\overline{\{ut\ \vert \ t\in \mathbb T\}}$ is compact, then
$u$ is called recurrent(in the sense of Birkhoff). The
corresponding motion $\pi (t, u)$ is then called a recurrent
motion. We identify the hull of the motion $\pi(t, u)$ with the
hull of the point $u$, i.e., $H(u)$.

 The $\omega-$limit set of $ u\in X$ is defined as
$\omega_u = \{v\; |\; \exists t_n \to +\infty \; \mbox{such that} \\
\;  \pi (t_n, u)\to v\}$.
 Likewise,  when $\mathbb T = \mathbb R $ or $\mathbb Z$,
 the $\alpha-$limit set of $u\in X$ is defined as
$\alpha_u = \{v \; |\; \exists t_n \to -\infty \; \mbox{such that}
\;  \pi (t_n, u)\to v\}$.

A point $u\in X$ is called stable in the sense of Poisson in the positive (negative)
direction if  $u\in \omega_u $ ($u\in \alpha_u $).
 A point $u\in X$ is called Poisson stable  if it is stable in the
 sense of Poisson in both the positive  and negative directions (in this case
 \ $\mathbb T = \mathbb R $ or $\mathbb Z$),
 i.e., $u\in \omega_u  \cap \alpha_u $.
 The corresponding motion $\pi (t, u)$ is
then called a Poisson stable motion.
 A recurrent motion is also
 a Poisson stable motion.
 Note that in the literature on topological  dynamics, a Poisson stable  point
is sometimes called a  recurrent point.

 A motion $\pi(t, u)$ is called pre--compact, if
the closure (under the topology of $X$) of its trajectory
 $\gamma(u)$, i.e., the hull $H(u):= {\mathcal Cl} \gamma(u)$, is compact.

 A set $E  \subset  X$  is
called invariant if the trajectory $\gamma (u)  \subset  E$
whenever $u \in E$.

 A set  $E  \subset  X$ is called minimal if it is
nonempty, closed and invariant, and it contains no proper subset
with these three properties.

Denote by $C(\mathbb R, X)$ the set of all continuous functions
$\phi: \mathbb R \to X$, equipped with the compact-open topology,
i.e., uniform convergence on compact subsets of $\mathbb R$. This
topology
 is metrizable as it can be generated by the complete metric
 $$
 d(\phi_1, \phi_2) = \sum_{n=1}^{\infty} \frac1{2^n}
 \frac{d_n(\phi_1, \phi_2)}{1+d_n(\phi_1, \phi_2)},
 $$
 where $
 d_n(\phi_1, \phi_2)=\max_{|t|\leq n}  \rho(\phi_1(t),\phi_2(t))$
  with $\rho$ the metric on $X$.

Let $\tau \in \mathbb R$ and $\phi_{\tau}$ be the
$\tau-$translation of function $\phi$, i.e., $\phi_{\tau} (t) :=
\phi(t+\tau)$ for all $t\in \mathbb R$. Define a mapping $\sigma:
\mathbb R \times C(\mathbb R, X) \to C(\mathbb R, X)$ by   $\sigma
(\tau, \phi):= \phi_{\tau}$. The triplet $( C(\mathbb R, X),
\mathbb R, \sigma)$ is a dynamical system and it is called a
dynamical system of translations (shifts) or the Bebutov's
dynamical system; see \cite{Sell, Scer1, Scer2}.
 Note that Fu and Duan  have shown that the Bebutov's
shift dynamical system is a   chaotic system \cite{FuDuan}.

A function $\phi \in C(\mathbb R, X)$ is   called  recurrent
(Poisson stable) if  it  is  recurrent (Poisson stable) under the
  Bebutov's dynamical system $( C(\mathbb R, X),\linebreak \mathbb R, \sigma)$.
The corresponding motion $\sigma(t, \phi)$ is also called
recurrent (Poisson stable) .

The set $H(\phi) := {\mathcal Cl} \{ \sigma(t, \phi) | t\in
\mathbb R\}$ is called the hull of the function $\phi$, where the
closure ${\mathcal Cl}$ is taken in the space $C(\mathbb R, X)$
with the compact-open topology.

\bigskip
We now discuss the difference between almost periodicity,
recurrence and Poisson stability. It is known (\cite{Sell, Si})
that
  an almost periodic point is recurrent, and a recurrent
point is Poisson stable. However, the converse    is not true. In
fact, let $\varphi \in C(\mathbb R, \mathbb R)$ be the function
defined by
$$ \varphi (t)= \frac{1}{2+\sin( t) +\sin (\pi t)}$$
for all $t\in \mathbb R$. It can be  verified that this function
is Poisson stable, but not recurrent, under the
  Bebutov's dynamical system, because it is not bounded.

An example of recurrent function which is not almost periodic can
be found, for example, in  Shcherbakov \cite{Sch72}.

In general, there exists the following relation between almost
periodicity (in the sense of Bohr) and recurrence (in the sense of
Birkhoff).
 Let $X $ be a complete metric space with metric $\rho$ and $(X, \mathbb R, \pi)$
be a dynamical system. A point $x\in X$ is almost periodic
  if and only if the following conditions are fulfilled: (1)
$x\in X$ is recurrent;  (2) a dynamical system $(H(x), \mathbb R,
\pi)$ is equicontinuous, i.e., for all $\varepsilon >0$ there
exists $\delta (\varepsilon) >0$ such that $\rho (x_1,x_2)<\delta$
implies $\rho (x_1 t, x_2 t) < \varepsilon$ for all $x_1,x_2 \in
H(x)$ and $t\in \mathbb R$. Note that the hull $H(x)=Cl\{xt : t
\in \mathbb R\}$.

Note that a function $\varphi \in C(\mathbb R, E)$ is almost
periodic if and only if the hull of the function is compact in the
uniform topology (with respect to $\sup$-norm).  However,  in the
present paper, we consider a compact-open topology for the
Bebutov's dynamical system.

\bigskip

Now we consider basic concepts in nonautonomous dynamical systems.

Let $(X,\mathbb R_{+},\pi)$ and $(\Omega ,\mathbb R ,\sigma)$ be
autonomous dynamical systems, and let $h: X \to \Omega $ be a
homomorphism from   $(X,\mathbb R _{+},\pi) $ onto $(\Omega
,\mathbb R ,\sigma)$.
 Then the triplet $\langle (X,\mathbb R _{+},\pi),(\Omega ,\mathbb
R ,\sigma),h\rangle $  is called   a {\em nonautonomous} dynamical
system (see \cite{Bro84},\cite{Ch02}). Since $(\Omega ,\mathbb R
,\sigma)$ is  usually a simpler  ({\em driving}) dynamical system,
the homomorphism $h$ may help understand the dynamical system
$(X,\mathbb R _{+},\pi) $.  This definition is slightly more
general than the definition used by some other authors.

Let $\Omega $ and $E$ be two metric spaces and $ (\Omega , \mathbb
R ,\sigma )$ be an autonomous dynamical system on $ \Omega.$
Consider a continuous mapping $ \varphi : \mathbb R_{+} \times E
\times \Omega  \to E $ satisfying the following conditions:
$$   \varphi (0,\cdot,\omega )=id_{E}, \ \ \varphi (t+\tau ,u,\omega)
=\varphi (t,\varphi (\tau,u,\omega ),\omega \tau) $$
for all $t, \tau \in \mathbb R ^{+}$, $\omega \in \Omega $ and $ u \in E $.
Such
 mapping $ \varphi $, or more explicitly, $\langle E, \varphi ,
(\Omega ,\mathbb  R ,\sigma ) \rangle $,  is called a cocycle on
the driving dynamical system $(\Omega, \mathbb R,\sigma ) $ with
fiber $E$; see \cite{Arn98, Sell}. This is a generalization of the
semi-flow property.

In the following example, we will see that the solution mapping
for the initial-value problem of a nonautonomous ordinary
differential equation usually defines a cocycle.

\begin{example}\label{ex3.1}
Let $E$ be a Banach space and $C(\mathbb R \times E,E)$ be a space
of all continuous functions $F: \mathbb R \times E \to E$ equipped
with the compact-open topology, i.e., uniform convergence on
compact subsets of $\mathbb R \times E$.
 Let us consider a
parameterized or nonautonomous differential equation
$$ \frac{du}{dt}+ Au = F(\sigma _t \omega ,u),   \ \  \omega \in \Omega, $$
on a Banach space $E$ with $\Omega$ $=$ $C(\mathbb{R}\times E,
E),$ where $\sigma _t\omega := \sigma (t,\omega)$ and
the linear operator $A$ is densely defined in
$E$ and such that the linear equation
\begin{eqnarray}\label{eq2.2*}
&& u'+ Au =0   \nonumber
\end{eqnarray}
generates a $c_0$-semigroup of linear bounded operators
$$
e^{-At}: E \to E, \  \varphi (t,x):=e^{-At}u .
$$
We will define
$\sigma_t$ $:$ $\Omega $ $\to$ $\Omega$
by $\sigma_t \omega (\cdot ,\cdot )$ $=$ $\omega (t+\cdot ,\cdot )$
for each $t$ $\in$ $\mathbb{R}$ and
interpret $\varphi(t,x,\omega)$ as mild solution of the initial value problem
\begin{equation} \label{eq3.1}
\frac{d }{dt} u(t) +Au(t) =  F(\sigma _t \omega ,u(t)), \quad u(0) = u.
\end{equation}
Under appropriate assumptions on $F$ $:$  $\Omega \times E$ $\to$
$E$,  or even $F$ $:$  $\mathbb{R}$ $\times E$ $\to$ $E $ with
$\omega (t)$ instead of $\sigma_t \omega $ in (\ref{eq3.1}), to
ensure forward  existence and uniqueness of the solution,
  $\varphi $ is a cocycle on
$( C(\mathbb R \times E, E), \mathbb R , \sigma) $ with
fiber $E$, where $( C(\mathbb R \times E, E), \mathbb R , \sigma) $
is a Bebutov's dynamical system
 which can be similarly defined as above;
 see, for example, \cite{Ch02}, \cite{Sch72},\cite{Sell}.
\end{example}

Note  that a cocycle induces a skew-product flow,
which naturally  defines a nonautonomous dynamical system.  In fact, let
$\varphi $ be a cocycle on $(\Omega ,\mathbb R ,\sigma )$ with the
fiber $E$. Then the mapping $\pi$ $:$
 $\mathbb R_{+}$  $\times E \times \Omega $  $\to$  $E \times \Omega $ defined by
$$
\pi(t,u,\omega) := (\varphi(t,u,\omega), \sigma _t \omega)
$$
for all $t$  $\in$  $\mathbb{R}_{+}$ and  $(u,\omega)$$\in$ $E\times
\Omega $ forms a semi-group
$\{\pi(t,\cdot,\cdot)\}_{t \in\mathbb{R}^{+}}$ of mappings of
$X:=\Omega \times E$
into itself,  thus a semi-dynamical system on the
 (expanded) state space $X$.   This semi-dynamical system is called a skew-product flow
(\cite{Sell}) and the triplet $\langle (X,\mathbb R
_{+},\pi),(\Omega ,\mathbb R ,\sigma),h\rangle $\ (where $h:=pr_2
: X \to \Omega )$
 is a nonautonomous dynamical system.

%%%%%%%%%%%%%%%%%%
%%%%%%%%%%%%%%%%%
 %%%%%%%%%%%%%%%%%
\section{Recurrent motions}

 We now  formulate the two-dimensional Navier-Stokes system
 as a nonautonomous dynamical system:
\begin{eqnarray}\label{eq2.9*}
& & u'+ \sum _{i=1}^{2}u_i\partial _i u =\nu \Delta u -\nabla p + F (t)  \\
& & div \ u =0, \ \ u|_{\partial D} =0 , \nonumber
\end{eqnarray}
where $ u= (u_1, u_2)$ is the velocity field, $F(x,y,t)=(F_1, F_2) $ is
the external forcing, and
$D $ is an open bounded fluid domain with smooth
boundary $ \partial D \in C^{2}$.

The functional setting of the problem is well known
\cite{Constantin,Temam}. We denote by $H$ and $V$ the closures of
the linear space $\{ u | u \in C_0^{\infty}(D)^{2} , \ div \ u=0
\} $ in $L^{2}(D)^{2}$ and $H_0^1(D)^{2}$, respectively. We also
denote by $P$ the corresponding orthogonal projection $P:
 L^2(D)^{2} \to H$. We futher set
$$
A:=-\nu P\Delta , \ B (u,v):= P(\sum_{i=1}^{2}u_i\partial_i v) .
$$

The Stokes operator $A$ is self-adjoint and positive with domain
$\mathcal D (A)$  dense in $H$. The inverse operator $A^{-1}$ is
compact. We define the Hilbert spaces $\mathcal D (A^{\alpha}), \
\alpha \in (0,1] $ as the domains of the powers of $A$ in the
standard way.  Note that $V:=\mathcal D (A^{1/2}),$ with norm
$\vert u \vert _{\mathcal D (A^{1/2})}=\vert  \nabla u\vert.$

Applying the orthogonal projection $P$, we rewrite (\ref{eq2.9*})
as an evolution equation in $H$
\begin{equation}\label{eq2.9**}
u'+Au +   B (u,u) = \mathcal F (t), \ \mathcal F (t):= P F (t) .
\end{equation}
We suppose that $\mathcal F\in C(\mathbb R ,H) $.
Let
  $(C(\mathbb R , H), \mathbb R , \sigma)$ be   the Bebutov's dynamical system;
see for example, \cite{Scer1},\cite{Scer2} and \cite{Sell}. We
denote $\Omega :=H(\mathcal F)=\overline{\{ \mathcal F_\tau | \
\tau \in \mathbb R \}}$,  where $ \mathcal F_\tau (t):=\mathcal F
(t+\tau) $ for all $t\in \mathbb R$ and the over bar   denotes the
closure in the compact-open topology in $C(\mathbb R , H)$. Note
that $\Omega$ is an invariant, closed subset of $C(\mathbb R ,
H)$.
 Moreover, $(\Omega ,\mathbb R, \sigma )$ is the Bebutov's dynamical
system of translations on $\Omega$.

Along with the equation (\ref{eq2.9**}),
 we consider the following family of equations
\begin{equation}\label{eq2.9***}
u'+Au +    B(u,u) = \tilde{\mathcal F} (t),
\end{equation}
where $ \tilde{\mathcal F}  \in \Omega=H(\mathcal F)$.
Let $ f:\Omega \to H $
be a mapping defined by
$$
 f(\omega)=f(\tilde{\mathcal F}):=
  \tilde{\mathcal F }(0)  ,
$$
where $\omega =  \tilde{\mathcal F} \in \Omega$. Thus, the
Navier-Stokes system can be written as a nonautonomous system
\begin{eqnarray}\label{NSE}
u_t+Au +B(u,u) &=& f(\omega t), \\
u(0) &=& u_0,
\end{eqnarray}
for $ \omega \in \Omega=H(\mathcal F)$ and $ u_0 \in H$. Here
$\omega t=\sigma_t(\omega)$. Note that $\Omega$ is compact minimal
if and only if $\mathcal F$ is recurrent. This nonautonomous
system is well-posed in $H$ for $t>0$.

 The solution $\varphi (t,u,\omega)$ of the nonautonomous Navier-Stokes
system (\ref{NSE})
  is called recurrent (Poisson stable, almost periodic,
quasi periodic),
if the point
$(u,\omega) \in H\times \Omega $ is a recurrent (Poisson stable,
almost periodic, quasi periodic)
point of
 the
skew-product dynamical system $(X,\mathbb R_{+},\pi)$, where
 $X=H\times \Omega$ and $\pi = (\varphi, \sigma)$.

As is known (in, for example, \cite{Sch85} and \cite{Lev-Zhi}), if
$\omega \in \Omega $ is a stationary ($\tau$-periodic, almost
periodic, quasi periodic, recurrent) point of dynamical system
$(\Omega ,\mathbb R ,\sigma)$ and $h:\Omega \to X $ is a
homomorphism of dynamical system $(\Omega ,\mathbb R ,\sigma)$
onto $(X,\mathbb R_{+},\pi)$, then the point $u=h(\omega)$ is a
stationary ($\tau$-periodic, almost periodic, quasi periodic,
recurrent) point of the system $(X,\mathbb R_{+},\pi)$.

It is known that (\cite{Vishik}) if the forcing function $\mathcal
F (t)$ is bounded with respect to time in $H$, such as in the case
of a recurrent function,  then every solution of the Navier-Stokes
system (\ref{eq2.9**}) is bounded on $\mathbb R_{+}$.

\begin{theorem}  ({\bf Recurrent motion})\\
If the forcing $\mathcal F(t)$ is recurrent in $C(\mathbb R, H)$,
then the 2D   Navier-Stokes system  (\ref{eq2.9**}) has at least
one recurrent solution in $C(\mathbb R, H)$.
\end{theorem}

\begin{proof}
The 2D   Navier-Stokes system (\ref{eq2.9**}) is  well-posed in
$H$. Moreover, the solution is in $H_0^1(D)$ (in fact, smooth)
after a short transient time due to regularizing effect; see for
example, \cite{Constantin, Foias2, Temam}. Since $H_0^1(D)$ is
compactly embedded in $H$, the solution is actually compact in $H$
after a short transient time.

Since the Navier-Stokes system (\ref{eq2.9**}) is nonautonomous,
its solution operator $\phi (t, u_0, g)$ does not define a usual
dynamical system or a semiflow in $H$.   However, we can define an
associated semiflow in an expanded phase space, i.e., a
skew-product flow  on $H \times H(\mathcal F)$:
$$
 \pi:  R_{+} \times H \times H(\mathcal F) \to  H \times H(\mathcal F),
$$ $$ \pi \ : \ (t, (u_0, \tilde{\mathcal F}))  \to  (\phi (t, u_0,
\tilde{\mathcal F}), \tilde{\mathcal F}_t).$$

Since $\mathcal F(t)$ is recurrent, $H(\mathcal F)$ is   compact
and minimal. Combining with the above compact solution $\phi (t,
u_0, \mathcal F)$ in $H$, we conclude that the skew-product flow
$\pi$  has
a pre-compact motion. Thus by a recurrence theorem due to Birkhoff
and Bebutov (see \cite{Sell}), the $\omega-$limit set
$\omega_{(u_0, \mathcal F)}$ contains a compact minimal set $M\subset
\omega_{(u_0,\mathcal F )}$. Let $(\tilde{u},\tilde{\mathcal F})\in M$.
Since $\mathcal F\in
H(\mathcal F)$ and $H(\mathcal F)$ is a compact minimal set,  we see that
$H(\tilde{\mathcal F})=H(\mathcal F)$. Thus there exists a sequence
$t_n \to +\infty$ such that $\tilde{\mathcal F}_{t_n}\to \mathcal F$.
Since the
sequence $\{\phi (t_n,\tilde{u},\tilde{\mathcal F})\}$ is pre-compact, then
without loss of generality we can suppose, that the sequence
$\{\phi (t_n,\tilde{u},\tilde{\mathcal F})\}$ is convergent. Denote by
$u:=\lim_{n\to \infty} \phi (t_n,\tilde{u},\tilde{\mathcal F})$, then the
point $(u,\mathcal F)\in M$ is recurrent and, consequently the
solution $\phi (t, u, \mathcal F) $ of equation (\ref{eq2.9**}) is
recurrent. This completes the proof.
 \end{proof}

\begin{remark}
This result is actually true for a more general evolution equation, as long as it generates     a cocycle and the
 solution  operator  is compact (or
asymptotically compact).

\end{remark}

\bigskip

We now assume that  the external forcing $\mathcal F(t) \in
C(\mathbb R, H)$ is
 Poisson stable.
We will show that the Navier-Stokes system (\ref{eq2.9**}) has at
least one Poisson stable  motion or solution.

We define $\mathfrak{N}_{\mathcal F}:= \{ \{t_n\} | \; \mathcal F
 t_n \to \mathcal F, t_n \to + \infty\}$.
 Let $\langle (X,\mathbb R _{+},\pi),(\Omega ,\mathbb R ,\sigma),\linebreak h\rangle $
be a  nonautonomous dynamical system. A point $u\in X$ is called
weakly regular (\cite{Sch72}), if for every $\{t_n\} \in
\mathfrak{N}_{y}$ with $y=h(u)$, the sequence $\{u t_n\}$ is
pre-compact in $X$. We need the following result.

\begin{theorem} ({\bf Shcherbakov \cite{Sch72}}) \\
Let $\langle (X,\mathbb R _{+},\pi),(\Omega ,\mathbb R
,\sigma),h\rangle $ be a  general nonautonomous dynamical system.
Let $y\in \Omega$ be Poisson stable under the dynamical system
$(\Omega ,\mathbb R ,\sigma) $ and there exists a weakly regular
solution $u_0$ of the abstract operator equation
$$
h(u)=y.
$$
Then this equation admits at least one Poisson stable solution in
$X$ under the dynamical system  $(X,\mathbb R _{+},\pi)$.
\end{theorem}

We have the following result.

\begin{theorem}({\bf Poisson stable  motion})\\
Let the forcing $\mathcal F(t)$ be Poisson stable  in $C(\mathbb
R, H)$. Assume that the Navier-Stokes system (\ref{eq2.9**})  has
one   bounded solution in $H_0^1(D)\cap H$ for $t>0$. Then the
Navier-Stokes system (\ref{eq2.9**}) has at least one Poisson
stable solution in $C(\mathbb R, H)$.
\end{theorem}

\begin{proof}
We again consider the skew-product flow  on $H \times H(\mathcal
F)$:
$$
 \pi:  R_{+} \times H \times H(\mathcal F) \to  H \times H(\mathcal F),
$$ $$ \pi (t, (u_0, \tilde{\mathcal F}))  \to  (\phi (t, u_0,
\tilde{\mathcal F}), \tilde{\mathcal F}_t). $$ Since
$\mathcal F(t)$ is  Poisson stable,   $\mathfrak{N}_{\mathcal F }
\not= \emptyset$. Let $\{t_n\}\in \mathfrak{N}_{\mathcal F} $.
Then the sequence $\{\mathcal F_{t_n}\}$ is convergent. By the
assumption,    the solution $ \phi (t, u_0,\mathcal F)$ is bounded
in $H^1_0$ on $R_{+}$. Consequently, due to the compact embedding
of $H^1_0$ in $H$,  the sequence $\{\pi (t_n, u_0, \mathcal F)) \}
$ is pre-compact in $H$. We argue that $(\phi (t_n, u_0, \mathcal
F), \mathcal F_{t_n})$ is pre-compact. Then by Theorem 3.2
we conclude that the
skew-product flow $\pi$ has at least one Poisson stable motion,
and hence the Navier-Stokes system (\ref{eq2.9**}) has at least
one
 Poisson stable solution.
This completes the proof.

 \end{proof}

\begin{remark}
This theorem also holds for a more general evolution
equation, as long as it generates a cocycle, admits at least
one bounded (on  $\mathbb R$) solution, and the
 solution  operator  is compact (or
asymptotically compact).

\end{remark}

\bigskip

{\bf Acknowledgment.} This work is partly supported by the Award
No. MM1-3016 of the Moldovan Research and Development Association
(MRDA) and U.S. Civilian Research \& Development Foundation for
the Independent States of the Former Soviet Union (CRDF),  and by
the NSF grant NSF-0139073. This paper was written while David
Cheban was visiting Illinois Institute of Thechnology (Departement
of Applied Mathematics) in Spring of 2002. He would like to thank
people there  for their very kind hospitality. We also thank the
referees for useful comments.

\medskip

Received June 2002; revised December 2002.

\medskip

\end{document}